\documentclass[12pt,a4paper]{amsart}
\usepackage[utf8]{inputenc}
\usepackage{amsmath, amssymb, amsthm}
\usepackage{enumerate}
\usepackage{lmodern}
\usepackage{xcolor}
\definecolor{darkgreen}{rgb}{0,0.5,0}
\usepackage{float}
\usepackage{graphicx}
\usepackage{xspace}
\usepackage{caption}
\usepackage{subcaption}
\usepackage{hyperref}
\hypersetup{
	colorlinks=true,
	citecolor=darkgreen,
	linkcolor=blue,
	filecolor=magenta,      
	urlcolor=cyan,
	pdftitle={Jacobian of a graph and graph automorphisms},
	pdfpagemode=FullScreen,
}

\newcommand{\red}[1]{\bgroup\color{red}{}#1\egroup}

\newtheorem{theorem}{Theorem}[section]
\newtheorem{problem}{Problem}
\newtheorem{lemma}[theorem]{Lemma}

\newtheorem{corollary}[theorem]{Corollary}

\theoremstyle{definition}
\newtheorem{definition}[theorem]{Definition}

\theoremstyle{remark}
\newtheorem{example}[theorem]{Example}

\numberwithin{equation}{section}
\newcommand{\aut}[1]{\ensuremath{\operatorname{Aut}(#1)}\xspace}
\newcommand{\mon}[1]{\ensuremath{\operatorname{Mon}(#1)}\xspace}
\newcommand{\homo}[1]{\ensuremath{\operatorname{H}_1(#1)}\xspace}
\newcommand{\jac}[1]{\ensuremath{\operatorname{Jac}(#1)}\xspace}
\newcommand{\sym}[1]{\ensuremath{\operatorname{Sym}(#1)}\xspace}
\newcommand{\fib}{\ensuremath{\operatorname{fib}}}

\makeatletter
\@namedef{subjclassname@2020}{%
  \textup{2020} Mathematics Subject Classification}
\makeatother

\title {The Jacobian of a graph and graph automorphisms}
\author[I. Est\'elyi]{Istv\'an Est\'elyi}
\address[I. Est\'elyi]{Faculty of Information Technology, University of Pannonia, Egyetem u. 10., 8200 Veszprém, Hungary}
\email[I. Est\'elyi]{estelyi.istvan@mik.uni-pannon.hu}
\author[J. Karab\' a\v s]{J\' an Karab\' a\v s}
\address[J. Karab\' a\v s]{Mathematical Institute, Slovak Academy of Sciences, Bansk\' a Bystrica, Slovakia}
\email[J. Karab\' a\v s]{karabas@savbb.sk}
\author[A. Mednykh]{Alexander Mednykh}
\address[A. Mednykh]{Sobolev Institute of Mathematics, Pr. Koptyuga 4, Novosibirsk, 630090, Russia}
\address[A. Mednykh]{Novosibirsk State University Pirogova 2, Novosibirsk, 630090, Russia}
\email[A. Mednykh]{smedn@mail.ru}
\author[R. Nedela]{Roman Nedela}
\address[R. Nedela]{Faculty of Applied Sciences, University of West Bohemia, Technická 8,
30100 Plzeň 3, Czech republic
}
\address[R.Nedela]{Mathematical Institute, Slovak Academy of Sciences, Bansk\' a Bystrica, Slovakia}
\email[R. Nedela]{nedela@kma.zcu.cz}

\begin{document}

\begin{abstract} In the present paper we investigate the faithfulness of certain linear representations of groups of automorphisms of a graph $X$ in the group of symmetries of the Jacobian of $X$. As a consequence	we show that if a $3$-edge-connected graph $X$ admits a nonabelian semiregular group of automorphims, then the Jacobian of $X$ cannot be cyclic. In particular, Cayley graphs of degree at least three arising from nonabelian groups have non-cyclic Jacobians. While the size of the Jacobian of $X$ is well-understood -- it is equal to the number of spanning trees of $X$ -- the combinatorial interpretation of the rank of Jacobian of a graph is unknown. Our paper presents a contribution in this direction. 
\end{abstract}


\maketitle

\section{Introduction}
\noindent{}The Jacobian of a finite graph is an important algebraic invariant  behaving nicely with respect to branched coverings of graphs~\cite{BN2007}. 
It is a certain finite Abelian group associated with the graph; we will introduce its formal definition in Section~\ref{sec:jacdef}. The notion  of the Jacobian of a graph, also known as the Picard group, the critical group, the dollar, or the sandpile group, was independently introduced by many authors (see e.g.~\cite{Bacher1997, baker2009, biggs1999, dhar1995, lorenzini2008, kotani2000}). It can be viewed as a discrete version of the Jacobian in the classical theory of Riemann surfaces. The Jacobian also admits natural interpretations in various areas of physics, coding theory, and financial mathematics. The fact that the size of the Jacobian of a connected graph $X$ is equal to the number of spanning trees is perhaps its most interesting property.

The present paper aims to investigate the relation between the symmetries of a graph $X$ and the symmetries of its Jacobian, $\jac{X}$. The main result is formulated in Theorem~\ref{thm:main}, where we prove that a semiregular group of symmetries of a connected, 3-edge-connected graph $X$  embeds into the automorphism group of $\jac{X}$. Corollary~\ref{cor:rank}
gives a sufficient condition for $\jac{X}$ to be acyclic, in particular, Jacobians of non-trivial Cayley graphs based on nonabelian groups are not cyclic (Corollary~\ref{cor:cay}).

The main idea of the proof of Theorem~\ref{thm:main} is to
derive two inequalities which  are contradictory. One of them is based on
properties of the Jacobians with respect to branched coverings of
graphs, for details see the paper \cite{BN2007} and Section~4. The other inequality follows from Lemma~\ref{lem:pfold} dealing with counting of spanning trees.
For technical reasons, we employ an extended model of
a graph, developed in Section~\ref{sec:defs}, allowing multiple edges, loops and even semiedges.

In \cite[Proposition~4.23]{baker2009}  it is proved that for an acyclic 2-connected graph $X$, the group $\aut{X}$ embeds into the group of symmetries of the homology group $\homo{X}$. On the other hand, there are  infinite families of graphs such that $\aut{X}$ does not embed into  $\homo{X}$. These families are completely determined in \cite[Theorem~7]{EKMN21}. Since $\jac{X}$ is a quotient of $\homo{X}$,  it is reasonable to include requirements on the connectivity of $X$.

\section{Graphs and graph coverings}
\label{sec:defs}
\noindent{}In this paper, we allow
graphs to have parallel edges, loops and semiedges. 
A simple undirected graph is defined in the standard way as a pair $X=(V,E)$, where $V$ is the set of vertices and $E$ is the set of edges. In this paper we prefer to define $X$ by means of a triple $(D;\sim,\lambda)$, where $D$
is a set of darts, $\sim$ is an equivalence relation on $D$ and $\lambda\in \sym{D}$ is an involutory permutation of $D$. An edge
$e=\{u,v\}$ gives rise to two darts $uv$ and $vu$, and for every edge $e\in E$ the dart-reversing involution $\lambda$ swaps the two underlying darts of $e$. The equivalence classes of $\sim$ determine the sets of darts $D_v$ incident to the vertex $v$, therefore one can view them as vertices of the graph.
The above-mentioned approach includes graphs with parallel edges, loops and possibly with semiedges.
Isolated vertices are not allowed in our model.   

Formally, a \emph{graph} $X$ is a triple $(D;\sim,\lambda)$, where $D$ is a finite, non-empty set of darts, $\sim$ is an equivalence on $D$, and $\lambda\in \sym{D}$ is an involution. The equivalence classes of $\sim$ will be called the \emph{vertices} of $X$, the set of all vertices will be denoted $V=V(X)$.
The orbits of $\lambda$ will be called the \emph{edges} of $X$. The projection $I\colon x\mapsto [x]_\sim$ from the set of darts onto the set of vertices will be called the \emph{incidence function}. For simplicity,
we often set
$\lambda(x)=x^{-1}$, for $x\in D$. 
An edge $e=\{x,x^{-1}\}$ is a \emph{semiedge}, if $x=x^{-1}$, the edge $e$ is
a \emph{loop} if $x^{-1}\neq x$, and $I(x^{-1})=I(x)$. Edges that are neither semiedges, nor loops, will be called \emph{ordinary edges}. Two
ordinary edges $e=\{x,x^{-1}\}$, $f=\{ y,y^{-1}\}$ will be called
\emph{parallel}, or multiple edges if $\{I(x),I(x^{-1})\}=\{I(y),I(y^{-1})\}$.
A graph $X$ will be called \emph{simple}, if every edge of $X$ is ordinary,
and $X$ has no parallel edges. The \emph{degree} of a vertex $v\in X$ is set to be $\operatorname{deg}(v)=|\{y\colon y\in D, v=I(y)\}|$. A \emph{walk} $W$ is a sequence
$x_1,x_2,\dots,x_{m}$ such that for every $i\in \{1,\ldots,m-1\}$
it holds that $I(x_i)=I(x_{i-1}^{-1})$. Graph $X=(D;\sim,\lambda)$ is \emph{connected} if for any two darts $x,y\in D$ there exists a walk $W = x_1,x_2,\dots,x_{m}$ such that $x=x_1$ and $y=x_m$. A walk $W$ is \emph{closed} if $I(x_1)=I(x_{m}^{-1})$. In particular, a single semiedge $x$ forms a trivial closed walk. A nontrivial closed walk $W$ is a (simple) \emph{cycle} if for any $i,j$, $i\neq j$ implies $I(x_i)\neq I(x_j)$. Note that the walks are by definition oriented. Further properties of simple graphs naturally extend to the present model of graphs and we will skip their definitions.

Let $X_1=(D_1;\sim_1,\lambda_1)$ and $X_2=(D_2;\sim_2, \lambda_2)$ be graphs. A mapping $\psi\colon D_1\to D_2$ is a \emph{graph homomorphism} 
if for any $x,y\in D_1$,
$x\sim_1 y$ implies $\psi(x)\sim_2\psi(y)$, and $\lambda_2(\psi(x))=\psi(\lambda_1(x))$.  In particular,
a permutation $f$ of the set of darts is an \emph{automorphism}
of a graph $X=(D; \sim,\lambda)$ if for every two darts $x,y\in D$,
$x\sim y$ implies $f(x)\sim f(y)$, and $\lambda(f(x))=f(\lambda(x))$. The homomorphism $\psi\colon X_1\to X_2$ is a \emph{covering}, if it is onto and its restriction to any class $[x]_{\sim_1}$ is a bijection. A covering
of connected graphs is by definition surjective.
The operators $I$ and $\lambda$ extends naturally to the set
of walks of $X_1=(D_1;\sim_1,\lambda_1)$. In particular, if $W=x_1,x_2,\dots ,x_m$ is a walk, then $I(W):=I(x_1)$ and $\lambda(W):=W^{-1}=x_m^{-1},x_{m-1}^{-1},\dots,x_1^{-1}$.
The set of preimages $\fib(z)=\psi^{-1}(z)$ for some $z\in D_2$ will be called the \emph{fibre} over $z$, similarly, by $\fib(v)$ we denote the \emph{vertex-fibre} over the vertex $v\in V(X_2)$. If $X_2$ is connected, then by the unique walk lifting property~\cite[Proposition~4.2]{MNS2000}, there exists $n$ such that $|\fib(z)|=n=|\fib(v)|$, 
for every dart $z\in D_2$ and for every vertex $v\in V(X_2)$.
Given a covering $\psi\colon X_1\to X_2$, the \emph{group of covering
	transformations} $\operatorname{CT}(\psi)\leq \aut{X_1}$ consists of those automorphisms $f$ satisfying $\psi\circ f = \psi$.
For each closed walk $W$ in $X_2$ originating at $v\in D_2$
we define the permutation $\beta_W$ of $\fib(v)$ by setting
$\beta_W(u)=I_1(\tilde W^{-1})$, where $\tilde W$ is the lift of a closed walk $W$
originating at $u$. The permutations $\beta_W$ acting
on $\fib(v)$ form a group called the \emph{monodromy group} $\mon{p}$.
If both $X_1$ and $X_2$ are connected, then the monodromy groups for different fibres
are isomorphic and the action of the covering transformation group on $\fib(u)$ and $\fib(v)$  are isomorphic for every $u,v\in V(X_2)$. The following two statements are well-known,
see \cite{MNS2000} for instance.

\begin{theorem}\label{thm:regcovers} Let $\psi\colon X_1\to X_2$ be an $n$-fold covering of connected
	graphs. Then
	\begin{enumerate}[{\rm (i)}]
		\item the monodromy group is transitive on each fibre,
		\item $\operatorname{CT}(\psi)$ is the centraliser of the monodromy action,
		\item the group of covering transformations is semiregular both on
		vertices and on darts, and
		\item $|\operatorname{CT}(\psi)|\leq n\leq |\mon{\psi}|$.
	\end{enumerate}
\end{theorem}

\begin{theorem}\label{thm:regcov} Let $\psi\colon X_1\to X_2$ be an $n$-fold covering of connected
	graphs. Then the following statements are equivalent
	\begin{enumerate}[{\rm (i)}]
		\item the monodromy group is regular on each fibre,
		\item $\operatorname{CT}(\psi)$ is regular on each fibre,
		\item $\operatorname{CT}(\psi)\cong \mon{\psi}$,
		\item $|\operatorname{CT}(\psi)|=n= |\mon{\psi}|$.
	\end{enumerate}
\end{theorem}
A covering satisfying any of the equivalent conditions in Theorem~\ref{thm:regcov} is called \emph{regular}. Note that if $\operatorname{CT}(\psi)\cong\mon{\psi}$ is an Abelian group, then $\operatorname{CT}(\psi)=\mon{\psi}$.

\begin{example}[Cayley graphs] Let $G=\langle S\rangle$
	be a group generated by a set $S=S^{-1}$ of non-trivial elements. Set $D=G\times S$, $\lambda(g,s)=(gs,s^{-1})$, and $(g,s)\sim (h,t)$
	when $g=h$. We have defined the Cayley graph $Cay(G;S)$ with respect
	to the generating set $S$. Let $B(S)=(S;\bar\sim,\bar \lambda)$,
	where $\bar\sim$ is the one-class equivalence and $\bar\lambda(s)=s^{-1}$. Clearly $B(S)$ is a one-vertex graph. Moreover, the projection
	$(g,s)\mapsto s$ is a regular covering  $Cay(G;S)\to B(S)$, with the group of covering transformations isomorphic to $G$. Since $S$ contains no trivial element, $Cay(G;S)$ has no loops. Further, since the action by left multiplication (on $G$) is semiregular, $Cay(G;S)$ contains no semiedges. Observe that the degree of every vertex of a Cayley graph $Cay(G;S)$ equals to $|S|$.
\end{example}

We have seen that each Cayley graph is a regular cover over a one-vertex graph, possibly with semiedges. It is natural to ask what is the class of graphs regularly covering one-vertex graphs. To answer the question one has to generalise the definition of a Cayley graph as follows. A \emph{Cayley multigraph} arising from a group $G$, $Cay(G,M)=(D;\sim,\lambda)$, is given by a multiset $M=M^{-1}$ of non-trivial elements of $G$, where $D=G\times M$, $\lambda(g,x)=(gx,x^{-1})$, and $(g,x)\sim (h,y)$
if and only if $g=h$. In particular, there are no loops and no semiedges in $Cay(G,M)$; however there may be parallel edges corresponding to the elements of $M$ of multiplicity greater than one.

The following lemma can be understood as a natural generalisation of the famous Sabidussi theorem~\cite{sabidussi}.
\begin{lemma}
	Let $X$ be a graph without loops and semiedges. Then $X$ regularly covers a one-vertex graph with the group of covering transformations $G$ if and only if $X$ is a Cayley multigraph arising from $G$.
\end{lemma}

The construction of a Cayley graph as a regular
cover over a one-vertex graph is generalised as
follows. Let $Y=(D;\sim,\lambda)$ be a connected graph and $T\subseteq Y$ be its spanning tree. Let $G$ be a finite group. A \emph{$T$-reduced
	voltage assignment} $\xi\colon D\to G$ is a mapping satisfying the following properties:
\begin{itemize}
	\item $\xi(x^{-1})=(\xi(x))^{-1}$,
	\item $\langle \xi(x),\ x\in D(Y)\rangle=G$,
	\item if $x\in D_T$, then $\xi(x)=1$.
\end{itemize}
The derived graph $X=Y^\xi=(\tilde D,\tilde\sim,\tilde\lambda)$ 
is defined by setting $\tilde D=G\times D$,
$(g,x)\tilde\sim (h,y)$ if and only if $x\sim y$
and $g=h$, $\tilde\lambda(g,x)=(g\xi(x),\lambda(x))$.

The following statement well-known, for graphs without semiedges \cite[Section~2.1]{gross2001}, is proved in \cite[Section 6]{MNS2000}.

\begin{theorem}\label{thm:va} Let $X$, $Y$ be connected graphs,
	$T\subseteq Y$ be a spanning tree,
	and $\xi\colon D(Y)\to G$ be a $T$-reduced voltage
	assignment. Then the projection $\psi\colon(g,x)\mapsto x$ is a regular covering $Y^\xi\to Y$ with the group of covering transformation isomorphic to $G$.
	
	Moreover, every regular covering 
	$X\to Y$ with the group of covering transformations isomorphic to $G$ can be described by means of a $T$-reduced
	voltage assignment $\xi\colon D(Y)\to G$.
\end{theorem}

Let $X$ be a connected graph. We denote the number of its spanning trees by $\tau(X)$.
The following lemma gives a lower bound for the number of spanning trees for a prime-fold cover
over a connected graph. It will be used
in the proof of Theorem~\ref{thm:main}.

\begin{lemma}\label{lem:pfold} Let $p$ be a prime and let $X$ be a simple, connected, $2$-edge-connected regular $p$-fold cover over a connected graph $Y$. Then $\tau(X)\geq p\cdot\tau(Y)$. If $X$ is $3$-edge-connected, then $\tau(X)>p\cdot\tau(Y)$.
\end{lemma}

\begin{proof} Let $\psi\colon X\to Y$ be a $p$-fold covering satisfying the assumptions. 
	By Theorem~\ref{thm:va} there is an associated $T$-reduced voltage assignment $\xi\colon D(Y)\to \mathbb{Z}_p$ such that the natural projection $\varphi\colon Y^\xi\to Y$ is equivalent to $\psi$.
	Hence, we may assume that $X=Y^\xi$ and 
	$\psi=\varphi$.
	Denote by $ F_T=\psi^{-1}(T)$.
	Clearly, $F_T$ is a spanning forest of $X$, consisting of $p$ isomorphic copies of $T$. If $p>2$, then since $X$ is simple, there are no semiedges in $Y$. Since $X$ is connected, there exists a cycle $C\subseteq Y$ (maybe a loop) which lifts to a cycle $\tilde C$ of length 
	$p\cdot|C|$.  
	Since $C$ lifts
	non-trivially, there exists a co-tree dart $x\in D(C)$ endowed with a non-trivial voltage. Let $e\in E(C)\setminus E(T)$ be the edge that contains $x$.
	Then $F_T+ \psi^{-1}(e)$ is an unicyclic spanning
	subgraph of $X$. Deleting any edge $g\in f^{-1}(e)$ from $F_T+ \psi^{-1}(e)$  we get a spanning tree $\tilde T$ of $X$. Removing different edges
	$g\in f^{-1}(e)$ we obtain a set $S(T,e)$ of $p$ spanning trees of $X$. We claim that 
	for different spanning trees $T_1\neq T_2$ of $Y$ we
	have $S(T_1,e_1)\cap S(T_2,e_2)=\emptyset$.
	If $\tilde T$ is a spanning tree constructed above,
	then it admits a unique decomposition of the 
	edge set into $|V(Y)|-1$ fibres over edges of $Y$, each of size $p$, and one incomplete fibre of size $p-1$. By definition, the complete fibres determine edges of a spanning tree $T$ of $Y$, and the incomplete fibre is a subset of $\psi^{-1}(e)$, where $e$ is a co-tree edge with a non-trivial voltage in the $T$-reduced assignment. Since the edge-decomposition of $\tilde T$ is unique, the spanning tree $T$ of $Y$ is uniquely determined as well.
	It follows that the sets $S(T,e)$, where $T$ ranges through the all spanning trees of $Y$ are pairwise disjoint, and we conclude
	$\tau(X)\geq p\cdot\tau(Y)$.
	
	Let $p=2$. If there exists a simple cycle in $Y$ which lifts nontrivially, we apply the same argument as above to prove that $\tau(X)\geq 2$. It may happen that there is no cycle in $Y$ which lifts non-trivially. Since $X$
	has no semiedges, every semiedge lifts to an ordinary edge.
	Since $X$
	is connected and $2$-edge-connected, $Y$ has at least two  semiedges $s_1$ and $s_2$ which lift
	to ordinary edges $e_1$ and $e_2$, respectively.
	Given spanning tree $T$ of $Y$, there are two associated spanning trees $F_T+e_1$ and $F_T+e_2$ of $X$. Since all spanning trees of $X$, constructed in this way are pairwise different, we have at least $2\cdot\tau(Y)$ spanning trees of $X$.
	
	Assume $X$ is $3$-edge-connected. If $p>2$, then given spanning tree $T$, there are at least two co-tree
	edges $e_1$, $e_2$ (maybe loops) assigned by a non-trivial voltage. If there was just one such edge, then $X$ would contain
	edge-cuts of size two. For a spanning tree $T$ of $Y$ we construct
	the sets $S(T,e_1)$ and $(T,e_2)$ employing
	the co-tree edges $e_1$ and $e_2$. Similarly, as above, we argue that the sets $S(T,e_i)$, $i=1,2$,
	and $T$ ranges through all spanning trees of $Y$,
	are pairwise disjoint. Hence, $\tau(X)\geq 2p\tau(Y)>p\tau(Y).$
	If $p=2$, and  there are two co-tree ordinary edges or loops $e_1$ and $e_2$, endowed with a non-trivial voltage,  we proceed as above.  Otherwise, since $X$ is $3$-edge-connected,
	either there are three semiedges $s_1$, $s_2$ and $s_3$ in $Y$, or
	there is a semiedge $s$ and a co-tree ordinary edge (or a loop) $e$ with a non-trivial voltage assignment.
	Since $X$ has no semiedges, every semiedge of $Y$
	lifts to an ordinary edge. In the first case,
	the sets $S(T,s_i)$, $i=1,2,3$, are singletons. In the second case,
	$S(T,e)$ is of cardinality two, and $S(T,s)$ contains one tree. Again, the sets are pairwise disjoint. It follows that $\tau(X)\geq 3\cdot\tau(Y)>2\cdot\tau(Y)$. 
\end{proof}

The following classical result on the edge connectivity of vertex-transitive graphs was proved independently by Mader and by Watkins.
\begin{theorem}[Mader~\cite{Mader70}, Watkins~\cite{Watkins70}] \label{Thm:VTkedgconn} The edge-connectivity of a vertex-transitive graph is equal to its valency.
\end{theorem}

\section{Jacobian of a graph}
\label{sec:jacdef}
\noindent{}The aim of this section is to introduce the Jacobian of a connected graph $X$ and to summarise some of its properties.  


\begin{definition} Let $X=(D;\sim,\lambda)$ be a connected graph.
	A mapping $\nu\colon D\to A$ into an Abelian group $A$ will be called
	an \emph{$A$-flow} if the following conditions are satisfied:
	\begin{description}
		\item[(FLW)] $\nu(x^{-1})=-\nu(x)$, for every $x\in D$ and
		\item[(GEN)]  $A$ is generated by $\{\nu(x)\ |\ x\in D\}$.
	\end{description}
	If an $A$-flow $\nu$ satisfies
	\begin{description}
		\item[(KLV)] $\sum_{y\sim x} \nu(y)=0$ for every dart $x\in D$ and 
		\item[(KLC)] $\sum_{i=0}^{m-1}\nu(x_i)=0$ for every oriented cycle $C=(x_0,x_1,\dots,x_{m-1})$,
	\end{description}
	it will be called a \emph{harmonic flow}.
\end{definition}

Note that the sets of equations (KLC) and (KLV) are the well-known Kirchhoff's laws for cycles and for vertices, respectively. In classical graph theory it is  usually required that an $A$-flow satisfies (KLV)~\cite[Chap.~6]{Diestel}. However, for our purposes the aforementioned definition is more appropriate.

\begin{definition} Given a connected graph $X$, the \emph{Jacobian $\jac{X}$}
	of $X$ is the maximal Abelian group $A$ such that $X$ admits
	a harmonic $A$-flow. A harmonic flow $D(X)\to \jac{X}$ will be called
	\emph{$J$-flow}.
\end{definition}

The Jacobian \jac{X} is not maximal just in ``numerical sense''.
If $A$ is an Abelian group such that there exists a harmonic $A$-flow $D(X)\to A$, then $A$ is a quotient of \jac{X}.

Denote by $\mathcal{A}(S)$ the free Abelian group generated by $S$.
Denote by $D^+$ a subset of $D$ containing from each edge 
$\{x,x^{-1}\}$ exactly one dart. In other words, $D^+$ is
a transversal of the set of edges of $X$. Observe that by (FLW)
a flow on $X$ is determined by its values on $D^+$. The following lemma gives a formal algebraic definition of a Jacobian.

\begin{lemma}\label{lem:universal} Let $X=(D;\sim,\lambda)$ be a connected graph. Then $\jac{X}\cong\mathcal{A}(D^+)/L$, where $L$ is the subgroup generated by the elements $\sum_{x\in C} x$ for every cycle $C$ of $X$ and by the elements $\sum_{y\sim z} y$, for every $z\in D$.
\end{lemma}

In what follows we list some properties of $\jac{X}$.

\begin{enumerate}[(P1)]
	\item $|\jac{X}|=\tau(X)$, in particular, $\jac{X}$ is a finite Abelian group~\cite{Bacher1997},
	\item Any Abelian group $A$, in which one can define a harmonic $A$-flow $D(X)\to A$, is a quotient of $\jac{X}$,
	\item $\jac{X}$ is a quotient of the homology group 
	$\homo{X, \mathbb{Z}}$,
	\item The rank of $\jac{X}$ is at most $v(X)-1$, where $v(X)$ is the number of vertices, and
	\item Finding a canonical decomposition  of $\jac{X}$ into cyclic factors is equivalent to computing the Smith normal form of
	the matrix of the homogeneous system of equations determined by $L$ over $\mathbb{Z}$~\cite{lorenzini2008}.
\end{enumerate}

Let us discuss the properties (P1), (P2) and (P3) which will be used in the next section. 

Property (P2) is a consequence of Lemma~\ref{lem:universal}.
Indeed, if $\nu\colon D(X)\to A$ is a harmonic $A$-flow, then by (GEN) the group $A$ can be viewed as a quotient of the free Abelian group 
$\mathcal{A}(D^+)$. Since $\nu$ is harmonic, (KLC) and (KLV) are satisfied.

Property (P3) is based on the following observation: Recall that $\homo{X, \mathbb{Z}}\cong \mathcal{A}(D^+_T)$, where  $D^+_T = D^+-D(T)$ with
respect to a spanning tree $T$ (co-tree darts). Since every tree has pendant darts
(darts incident to vertices of degree one), if the values of a harmonic flow on $D^+_T$ are prescribed, then using
(KLV) one can extend the flow to the pendant darts in $D^+\setminus D^+_T$. Repeating the argument, the flow extends to all the darts in $D^+$. Hence, $\jac{X}$ is a quotient of $\homo{X, \mathbb{Z}}$. Each semiedge contributes to the homology group of $X$ by a $\mathbb{Z}_2$-factor and these factors are independent~\cite{MNS2000}. In particular, the semiedges behave as independent cycles and therefore by (KLC) they do not contribute to \jac{X}.

Property (P1) is a reformulation of the well-known Kirchhoff matrix-tree theorem.
To see this, one has to relate our definition of the Jacobian to the standard one, see \cite{BN2007}, using the concept of divisors. A \emph{divisor} of $X$ is an integer-valued function defined on the set of vertices of $V$. Denote by $\operatorname{Div}(X)$ the set of divisors of $X$. Observe that the sum of two divisors is again a divisor; therefore the divisors of $X$ form a free Abelian group of rank $v(X)$.
A divisor $f\colon V(X)\to\mathbb{Z}$ determines a flow
by setting $\nu_f(x)=f(I(x^{-1}))-f(I(x))$. The flow $\nu_f$ satisfies
(KLC) for all closed walks of $X$. Vice-versa, any flow $\nu$ in an Abelian group $A$ satisfying (KLC) determines a function $g\colon V\to A$ such that $\nu=\nu_g$. Note that $g$ is not uniquely determined
by $\nu$. To make this correspondence unique, we need to fix a
value of $g$ at some vertex (or to introduce any other linear relation on the set of vertices). It follows that the set of flows in a group $A$ satisfying (KLC) is in correspondence with the set of functions $V(X)\to A$ satisfying an extra linear relation. A good question to ask
is which functions on vertices correspond to flows satisfying both
(KLC) and (KLV). With a little effort one can establish the following correspondence.

\begin{lemma} Let $X=(D;\sim,\lambda)$ be a graph without loops and semiedges and $\Delta$ be its Laplacian. Let $f\colon V\to A$ be a function. Then
	a flow $\nu_f\colon D(X)\to A$ is a harmonic $A$-flow if and only if the function $f$ satisfies $\Delta f= \vec 0$. 
\end{lemma}

It follows that $\jac{X}\cong \operatorname{Div}_0/\Delta(\operatorname{Div})$, where $\operatorname{Div}_0$ is
the set of divisors of degree $0$. Recall that a divisor $f\colon V(X)\to A$ is of degree $0$ if $\sum_{v\in V(X)} f(v) = 0$. With this correspondence in mind, the proof of (P1) can be found in \cite[p.~769]{BN2007}. The equivalence of the various definitions of Jacobians, including the one used throughout the paper can be found in the appendix of~\cite{LLQZ20}.

\section{Representation of $\aut{X}$ in $\aut{\jac{X}}$}
\noindent{}Assume
we have a harmonic $J$-flow $\xi\colon D\to \jac{X}$. Then for every $f\in\aut{X}$ we can
define a transformation $\xi\mapsto \xi_f$ by setting $\xi(x)=\xi(f(x))$. Since $f$ takes vertices onto vertices and cycles onto cycles, $\xi_f$ is a $J$-flow. Since $\xi$ generates $\jac{X}$,
every element $a$ of the Jacobian can be written as $a=\sum_{x\in D}c_x\xi(x)$, where $c_x\in \mathbb{Z}$ for each $x\in D$. We define a mapping $\Theta\colon\aut{X}\to\aut{\jac{X}}$
taking $f\mapsto f^*$, $f\in\aut{X}$, where $f^*$ is given by setting 

$$f^*(a)=f^*(\sum_{x\in D}c_x\xi(x)):=\sum_{x\in D}c_x\xi(f(x)).$$

\begin{lemma}With the above notation, $f^*$ is an automorphism of $\jac{X}$ and $\Theta\colon\aut{X}\to\aut{\jac{X}}$ is a group homomorphism.
\end{lemma}

\begin{proof} The restriction of $f^*$ onto the generating
	set $S=\{\xi(x)\colon x\in D\}$ of $\jac{X}$ permutes the elements of $S$. Moreover, $f\in\aut{X}$ satisfies the following
	\begin{itemize}
		\item $\xi f(x^{-1})=\xi((f(x))^{-1})=-\xi f(x)$, for every $x\in D$,
		\item the set $D_v:=\{ x\in D \colon I(x)=v\}$ of darts based at $v\in V$ is mapped by $f$ onto $D_{f(v)}$, hence 
		$$f^*\left(\sum_{x\in D_v}\xi(x)\right)=\sum_{f(x)\in D_{f(v)}}\xi(f(x))=\sum_{z\in D_{f(v)}}\xi(z),\text{ and}$$
		\item an oriented cycle $C$ is mapped to an oriented cycle $f(C)$, hence
		$$f^*\left(\sum_{x\in C}\xi(x)\right)=\sum_{f(x)\in f(C)}\xi(f(x))=\sum_{z\in f(C)}\xi(z).$$
	\end{itemize}
	It follows that a relation in the presentation of $\jac{X}$ is mapped onto a relation. Consequently, the mapping $\xi(x)\mapsto \xi(f(x))$ extends
	to a group automorphism.
	
	Now we verify that $\Theta$ is a group homomorphism. Let $f,g\in\aut{X}$ be automorphisms. Then
	\begin{gather*}
		\Theta(f\circ g)(y)=(f\circ g)^*(y)=
		\sum_{x\in D}c_x\xi((f\circ g)(x))=\alpha^*\sum_{x\in D}c_x\xi(g(x))\\
		=f^*(g^*(y))=\Theta(f)(\Theta(g)(y)).\qedhere
	\end{gather*}
\end{proof}
The central problem we are interested in reads as follows: For which subgroups $G\leq\aut{X}$ is the image $\Theta(G)$ isomorphic to $G$?
Equivalently, our aim is to investigate under which conditions the restriction 
$\Theta|_G$ is a monomorphism. By the first isomorphism theorem,
a subgroup $G\leq \aut{X}$ embeds into $\aut{\jac{X}}$
if and only if the kernel $\ker(\Theta|_G)$ is trivial. 
Observe that $f\in \ker(\Theta|_G)$
if $\xi(x)=\xi(f(x))$ for every $x\in D$, i.e., if and only if $\xi$ takes constant values on the dart-orbits of the group $\langle f\rangle$, generated by $f$. An automorphism $f$ will be called \emph{$\xi$-invariant} if $\xi(x)=\xi(f(x))$, for every dart $x\in D$.

The following lemma deals with properties of $\xi$-invariant automorphisms. For a flow $\xi\colon D(X)\to A$ the \emph{local group} $A^\xi\leq G$ is the subgroup of $A$ generated 
by $\{\xi(C)\colon $C$ \text{\ an\ oriented\ cycle\ of $X$}\}$, where 
$\xi(C)=\sum_{i=0}^{m-1}\xi(x_i)$, for an oriented cycle $C=(x_0,x_1,\dots,x_{m-1})$.
If $\xi$ is a harmonic flow, then the local group is trivial by (KLC).

\begin{lemma}\label{lem:cycliccov} Let $X$ be a connected graph and let $\xi\colon D\to A$ be a harmonic flow on $X$. Let $f\in \aut{X}$ be a $\xi$-invariant automorphism which is semiregular both on darts and vertices. Let $\varphi\colon X\to X/\langle f\rangle$ be the canonical projection  $x\mapsto [x]_f$, for $x\in D(X)$. 
	Then $\bar \xi\colon \bar D\to A$ defined by $\bar\xi([x]_f):=\xi(x)$ is a flow satisfying {\rm (KLV)} on $Y:=X/\langle f\rangle$ and the local group $A^{\bar\xi}$ is an epimorphic image of $\langle f\rangle$. 
	
\end{lemma}

\begin{proof} 
	Let the order of $f$ be $n$.
	By definition, $\bar\xi$ is a flow satisfying (KLV).
	Since the projection $X\to Y$ is a regular covering, the vertices
	in the fibre over $v\in Y$ can be indexed as $v_0,v_1,v_2,\dots,v_{n-1}$, where $v_i:=f^i(v_0)$.
	By Theorem~\ref{thm:regcov}, the action of the covering transformation group $\langle f\rangle$ coincides with the monodromy action of the group of $v$-based closed walks $\pi(v,Y)$ on $\fib_v=\{v_0,v_1,v_2,\dots,v_{n-1}\}$.
	More precisely, there is an epimorphism $\Phi\colon \pi(v,Y)\to\langle f\rangle$ with 
	$\ker(\Phi)= \varphi(\pi(v_0,X))$. We show that the assignment: $\Psi\colon \Phi(W)\mapsto\bar\xi(W)$, $W\in\pi(v,Y)$, is an epimorphism from $\langle f\rangle$ to the local group
	$A^{\bar\xi}$.
	
	First we show that $\Psi$ is well-defined. Suppose that 
	$\Phi(W_1)=\Phi(W_2)$. Then
	$1=W_1W_2^{-1}$ and $\Phi(W_1W_2^{-1})=0$ implying
	that there are lifts $\tilde W_1$ and $\tilde W_2$ such that $\tilde W_1\tilde W_2^{-1}\in \pi(v_0,X).$ 
	Since $\xi$ is a harmonic flow, we have $\xi(\tilde W_1\tilde W_2^{-1})=0$, and consequently, $\bar\xi(W_1W_2^{-1})=0$. Thus, $\bar\xi(W_1)=\bar\xi(W_2)$.
	Hence, $\Psi(\Phi(W_1))=\bar\xi(W_1)=\bar\xi(W_2)=\Psi(\Phi(W_2))$.
	
	Secondly, we show that it is a homomorphism. Let $g_1,g_2\in\langle f\rangle$, where $g_1=\Phi(W_1)$ and
	$g_2=\Phi(W_2)$ for some $W_1,W_2\in \pi(v,Y)$. Then 
	\begin{align*}
		\Psi(g_1g_2)&=\Psi(\Phi(W_1)\Phi(W_2))=
		\bar\xi(W_1W_2)=\Psi(\Phi(W_1W_2))\\
		&=\bar\xi(W_1)\bar\xi(W_2)=\Psi(\Phi(W_1))\Psi(\Phi(W_2))=\Psi(g_1)\Psi(g_2).
	\end{align*}
	By definition, $\Psi$ is an epimorphism.
\end{proof}

Now we are ready to prove the main result of the paper establishing that any semiregular group of automorphisms of a graph  acts faithfully on its Jacobian. 

\begin{theorem}\label{thm:main} Let $X$ be simple, connected, $3$-edge-connected graph. Let $G\leq\aut{X}$ be a subgroup acting semiregularly both on darts and on vertices of $X$. Then the subgroup $\Theta(G)$ of $\aut{\jac{X}}$ is isomorphic to $G$.
\end{theorem}

\begin{proof} Let $\xi$ be a harmonic $J$-flow in $\jac{X}$.
	Our aim is to prove that  the kernel of the homomorphism 
	$\Theta|_G\colon G\to \aut{\jac{X}}$ is trivial. Suppose to the contrary that there exists $f\in G$, $f\neq \operatorname{id}$
	such that $f\in \ker(\Theta|_G)$. It follows that $\xi(f(x))=\xi(x)$
	for every dart $x\in D(X)$, in other words $f$ is $\xi$-invariant. Clearly, taking a proper power of $f$, we obtain an element of $\ker(\Theta|_G)$
	of prime order $p$. Hence, we may assume that
	$f$ has order $p$.
	Denote $\bar x=[x]_f$ 
	the orbit of $x\in D(X)$ under $\langle f\rangle$. 
	Since the action of $G$ is semiregular on darts and vertices, 
	the natural projection $x\mapsto [x]_f=\bar x$ is
	a regular cyclic covering $\gamma\colon X\to Y$, $Y:=X/\langle f\rangle$.
	By Lemma~\ref{lem:cycliccov}, $\gamma$ induces 
	a flow $\bar\xi\colon D(Y)\to \jac{X}$ on $Y$ satisfying (KLV). Moreover, the corresponding local group $(\jac{X}/\langle f\rangle)^{\bar{\xi}}$
	is cyclic of order dividing the order of $f$,
	that means it is of order $p$.
	The flow  $\bar\xi$  may not satisfy  (KLC). However,
	we can easily understand the ``defect'' of $\bar\xi$.
	For every oriented simple cycle $\bar C$ of $Y=X/\langle f\rangle$ we set
	$df(\bar C):=\sum_{\bar x\in D(\bar C)}\bar\xi(\bar x)$. The total
	defect of $\bar\xi$ is the subgroup 
	$$K:=df(\bar\xi):=\langle  df(\bar C)|\  \bar C\in {\mathcal C}(Y)\rangle,$$
	where  ${\mathcal C}(Y)$ is the set of oriented cycles in $Y$. Now we can formulate several claims about $K$:
	
	\begin{itemize}
		\item $K\leq \jac{X}$,
		\item the flow $\xi^*$ in the factor group $\jac{X}/K$, 
		defined by $\xi^*(\bar x):=\bar\xi(\bar x)K$, satisfies both Kirchhoff laws, and
		\item $K$ is a cyclic group whose order divides $|f|$.
	\end{itemize}
	
	The first claim holds true, since every generator  $df(\bar C)$ of $K$ belongs to $\jac{X}$. For every
	cycle $\bar C$, we have $\sum_{\bar x\in \bar C}\xi^*(\bar x)=\sum_{\bar x\in \bar C}\bar\xi(\bar x)K=df(\bar C)K=K$. Thus the second claim holds.  
	By definition, $K$ is a local group with respect to $\bar{\xi}$.
	By Lemma~\ref{lem:cycliccov}, $|K|$ divides $|f|=p$.
	Since $\xi^*$ is a harmonic flow in $\jac{X}/K$, 
	the group $\jac{X}/K$ is an epimorphic image of $\jac{Y}$. It follows that 
	$$\tau(Y)=\tau(X/\langle f\rangle)=|\jac{Y}|\geq |\jac{X}/K|=
	|\jac{X}|/|K|\geq\tau(X)/p.$$ 
	Thus, $\tau(X)\leq p\tau(Y)$. 
	However, Lemma~\ref{lem:pfold} gives $\tau(X)> p\tau(Y)$, a contradiction.
\end{proof}

\begin{example} Let $X$ be a graph with one central vertex and two outer vertices
	with three parallel edges between the central vertex and each outer vertex.
	Then $|\aut{X}|=72$ and $\jac{X}\cong {\mathbb Z}_3\times{\mathbb Z}_3$.
	It follows that $\aut{\jac{X}}$ is the general linear group $\operatorname{GL}_2(3)$  of order $48$. 
	Hence, $\aut{X}$ does not embed into $\aut{\jac{X}}$.
	The example can be easily generalised to an infinite family.
	Hence, the assumption of semiregularity of $G$ in Theorem~\ref{thm:main}
	is essential. 
\end{example}

\begin{corollary}\label{cor:rank} If a connected, $3$-edge-connected graph $X$
	admits a nonabelian semiregular group of automorphisms, then the rank of the Jacobian $\jac{X}$ is at least $2$.
\end{corollary}

\begin{proof} Assume, to the contrary, that $\jac{X}$ is cyclic.
	Then $\aut{\jac{X}}$ is Abelian. However, by Theorem~\ref{thm:main}
	the group $\aut{\jac{X}}$ contains a nonabelian subgroup, a contradiction.
\end{proof}

\begin{corollary}\label{cor:cay} Let $X$ be a Cayley graph arising from a nonabelian group. If the degree of $X$ is at least three, then $\jac{X}$ is not cyclic.
\end{corollary}
\begin{proof} By Theorem~\ref{Thm:VTkedgconn} the graph $X$ is
	$3$-edge-connected. Now the statement follows from
	Corollary~\ref{cor:rank}.
\end{proof}

It is natural to ask whether  Theorem~\ref{thm:main} can be generalised. In particular, the following open problem.
\begin{problem}
	Is there a $3$-connected  graph $X$  such that \aut{X} does not embed into \aut{\jac{X}}?
\end{problem}

Corollaries~\ref{cor:rank} and \ref{cor:cay} suggest that the structure of the automorphism group \aut{X} influences the rank of \jac{X}. In general we have the following general problem.
\begin{problem}
	Which properties of simple connected graphs are related to the rank of its Jacobian? 
\end{problem}

Computer-aided experiments with small graphs suggest that most graphs have cyclic Jacobians. The following problem is of interest.

\begin{problem}
	Characterise simple $2$-connected graphs with acyclic Jacobians.
\end{problem}
A sufficient condition for the group \jac{X} to be cyclic has been obtained in terms of the characteristic polynomial of the Laplacian matrix of $X$
in~\cite[Lemma~2.13]{lorenzini2008}.


\section*{Acknowledgements}
\noindent{}The authors wish to express their sincere gratitude to the anonymous referees for many valuable comments leading to significant improvements of the paper. 

The first, second and the fourth authors were supported by the grant GACR 20-15576S. The second and fourth authors were supported by the grant No. APVV-19-0308 of Slovak Research and Development Agency and grant VEGA 02/0078/20 of Slovak Ministry of Education. The third author was supported by Sobolev Institute of Mathematics under agreement No. FWNF-2022-0005 with the Ministry of Science and Higher Education of the Russian Federation.

\bibliographystyle{abbrv}
\bibliographystyle{abbrv}

\end{document}